\documentstyle[leqno]{book}
\def\mA{\mbox{\bf A}}
\def\mB{\mbox{\bf B}}
\def\mR{\mbox{\bf R}}

\begin{document}

\thispagestyle{plain}

\noindent {\small\sc Univ. Beograd. Publ. Elektrotehn. Fak.}

\noindent {\scriptsize Ser. Mat. 9 (1998), 29--33}

\vspace*{10.00 mm}

\centerline{\Large \bf SOME COMBINATORIAL ASPECTS}

\medskip
\centerline{\Large \bf OF DIFFERENTIAL OPERATION}

\medskip
\centerline{\Large \bf COMPOSITION ON THE SPACE
$\mbox{\Large \bf R}^{\mbox{\small \bf n}}$}
\footnotetext{1991 Mathematics Subject Classification: 26B12, 58A10}

\vspace*{4.00 mm}

\centerline{\large \it Branko J. Male\v sevi\' c}

\vspace*{4.00 mm}

\begin{center}
\parbox{25.0cc}{\scriptsize \bf In this paper we present a recurrent
relation for counting meaningful compositions of the higher-order differential
operations on the space $\mbox{\scriptsize \bf R}^{\mbox{\tiny \bf n}}$
(n=3,4,...) and extract the non-trivial compositions of order higher than two.}
\end{center}

\medskip

\medskip
\centerline{\bf 1. DIFFERENTIAL FORMS AND OPERATIONS
ON THE SPACE ${\mbox{\bf R}}^{\mbox{\footnotesize \bf 3}}$}

\bigskip
It is well known that the first-order differential operations grad, curl
and div on the space $\mR^3$ can be introduced using the operator of the
exterior~differentiation~$d$ of differential forms [{\bf 1}]:
{\small $$
\mbox{\normalsize $\Omega$}^{0}(\mR^{3})
\stackrel{d}{\longrightarrow}
\mbox{\normalsize $\Omega$}^{1}(\mR^{3})
\stackrel{d}{\longrightarrow}
\mbox{\normalsize $\Omega$}^{2}(\mR^{3})
\stackrel{d}{\longrightarrow}
\mbox{\normalsize $\Omega$}^{3}(\mR^{3}),
$$}
\hspace*{-2.5 mm} where $\Omega^{i}(\mR^{3})$ is the space of differential
forms of degree $i = 0, 1, 2, 3$ on the space $\mR^3$ over the ring of functions
{\small $\mbox{\normalsize \mA} = \{ f: \mR^3 \rightarrow \mR \, | \, f \in
C^{\infty}(\mR^3) \}$}. In the consideration, which follows, we give definitions
of the first-order differential operations.

\smallskip
Let us notice that one-dimensional spaces $\Omega^{0}(\mR^3)$ and
$\Omega^{3}(\mR^3)$ are isomorphic to $\mA$ and let  $\varphi_{0}:$
{\small $\Omega^{0}(\mR^3) \rightarrow \mA$}, $\varphi_{3}:$ {\small
$\Omega^{3}(\mR^3) \rightarrow \mA$} be the corresponding isomorphisms.
Next, the set of vector functions \mbox{\small $\mbox{\normalsize \mB}
= \{ \mbox{\boldmath $f$}\!=\!(f_1, f_2, f_3) : \mR^3 \rightarrow \mR^3 \, | \,
f_1,f_2,f_3 \in C^{\infty}(\mR^3) \}$}, over the ring $\mA$, is
three-dimensional. It is isomorphic to $\Omega^{1}(\mR^3)$~and
$\Omega^{2}(\mR^3)$. Let $\varphi_{1}\!\::\!\:$ {\small $\Omega^{1}(\mR^3)
\rightarrow \mB$}, $\varphi_{2}\!\::\!\:$ {\small $\Omega^{2}(\mR^3)
\rightarrow \mB$} be the corresponding isomorphisms. In that case,
the compositions $\varphi_{0}^{-1}\!\:\circ \varphi_{3}:$ \mbox{\small
$\Omega^{3}(\mR^3) \rightarrow \Omega^{0}(\mR^3)$} and $\varphi_{1}^{-1}\!\:
\circ \varphi_2:$ \mbox{\small $\Omega^{2}(\mR^3) \rightarrow \Omega^{1}(\mR^3)$}
are isomorphisms of the corresponding spaces of differential forms. The
first-order differential operations are defined via the operator of the
exterior differentiation $d$ of differential forms in the following form:
$$
\nabla_{1}\!\:=\!\:\varphi_{1}\!\:\circ\!\:d\!\:\circ\!\:\varphi_{0}^{-1}:
\mA \rightarrow \mB, \;\;
\nabla_{2}\!\:=\!\:\varphi_{2}\!\:\circ\!\:d\!\:\circ\!\:\varphi_{1}^{-1}:
\mB \rightarrow \mB, \;\;
\nabla_3\!\:=\!\:\varphi_3\!\:\circ\!\:d\!\:\circ\!\:\varphi_{2}^{-1}:
\mB \rightarrow \mA.
$$
Therefore we obtain explicit expressions for the first order differential
operations  $\nabla_{1}$, $\nabla_{2}$, $\nabla_3$ on the space $\mR^3$
in the following form:

{\small
\smallskip                                                                
\noindent {\normalsize (1)}
$\quad\displaystyle \mbox{\normalsize grad}\,\mbox{\normalsize $f$}=
\mbox{\normalsize $\nabla_1$} \mbox{\normalsize $f$}=
\frac{\partial f}{\partial x_1}\,\mbox{\normalsize \boldmath $e_1$}+
\frac{\partial f}{\partial x_2}\,\mbox{\normalsize \boldmath $e_2$}+
\frac{\partial f}{\partial x_3}\,\mbox{\normalsize \boldmath $e_3$}:
\mA \rightarrow \mB,$

\smallskip                                                               
\noindent {\normalsize (2)}
$\quad\displaystyle \mbox{\normalsize curl}\,\mbox{\normalsize \boldmath $f$}=
\mbox{\normalsize $\nabla_2$} \mbox{\normalsize \boldmath $f$}=
\left(\frac{\partial f_3}{\partial x_2}
\!-\!
\frac{\partial f_2}{\partial x_3}\right)\mbox{\normalsize \boldmath $e_1$}+
\left(\frac{\partial f_1}{\partial x_3}
\!-\!
\frac{\partial f_3}{\partial x_1}\right)\mbox{\normalsize \boldmath $e_2$}+
\left(\frac{\partial f_2}{\partial x_1}
\!-\!
\frac{\partial f_1}{\partial x_2}\right)\mbox{\normalsize \boldmath $e_3$}:
\mB \rightarrow \mB,$

\smallskip                                                               
\noindent {\normalsize (3)}
$\quad\displaystyle \mbox{\normalsize div}\,\mbox{\normalsize \boldmath $f$}=
\mbox{\normalsize $\nabla_3$} \mbox{\normalsize \boldmath $f$}=
\frac{\partial f_1}{\partial x_1}+
\frac{\partial f_2}{\partial x_2}+
\frac{\partial f_3}{\partial x_3}:
\mB \rightarrow \mA.$}

\break

     Let us count meaningful compositions of differential operations $\nabla_1,
\nabla_2, \nabla_3$. Consider the set of functions $\Theta = \{ \nabla_1,
\nabla_2, \nabla_3 \}$. Let us define a binary relation $\rho$~{\it "to be
in composition"} with $\nabla_i \rho \nabla_j = \top$ iff the composition
$\nabla_{j} \circ \nabla_{i}$ is meaningful $(\nabla_i, \nabla_j \in
\Theta)$. The {\sc Cayley}'s table of this relation reads:
{\small $$                                                              
\begin{array}{c|rrr}
\rho       & \nabla_{1} & \nabla_{2} & \nabla_{3} \\ \hline
\nabla_{1} & \bot       & \top       & \top       \\
\nabla_{2} & \bot       & \top       & \top       \\
\nabla_{3} & \top       & \bot       & \bot
\end{array}_{\mbox{\normalsize .}}
\leqno \mbox{\normalsize (4)}
$$}
\hspace*{-2.5 mm} We form the graph of relation $\rho$ as follows. If
$\nabla_i \rho \nabla_j = \top$ then we put the node $\nabla_j$ under
the node $\nabla_i$. Let us mark $\nabla_0$ as nowhere-defined function
$\vartheta$, with domain and range being the empty set [{\bf 2}]. We shall
consider $\nabla_0 \rho \nabla_i = \top$ $(i=1,2,3)$. For the set of
functions $\Theta \cup \{\nabla_{0}\}$ our graph is the tree with
the root in the node $\nabla_0$.

\vspace*{-4 mm}                                                             
\setlength{\unitlength}{0.17 cc}                                            
\begin{picture}(150,50)(0,0)                                                
\thicklines                                                                 
\put(60,40){\circle*{0.8}}                                                  
\put(61,41){\scriptsize$\nabla_{0}$}                                        
\put(130,41){\small$f(0)=\;1$}                                              
\put(20,30){\line(4,1){40}}                                                 
\put(20,30){\circle*{0.8}}                                                  
\put(17,31){\scriptsize$\nabla_{1}$}                                        
\put(60,30){\line(0,1){10}}                                                 
\put(60,30){\circle*{0.8}}                                                  
\put(61,31){\scriptsize$\nabla_{2}$}                                        
\put(100,30){\line(-4,1){40}}                                               
\put(100,30){\circle*{0.8}}                                                 
\put(101,31){\scriptsize$\nabla_3$}                                         
\put(130,31){\small$f(1)=\;3$}                                              
\thinlines                                                                  
\put(10,20){\line(1,1){10}}                                                 
\put(10,20){\circle*{0.8}}                                                  
\put(6,21){\scriptsize$\nabla_{2}$}                                         
\thicklines                                                                 
\put(30,20){\line(-1,1){10}}                                                
\put(30,20){\circle*{0.8}}                                                  
\put(31,21){\scriptsize$\nabla_3$}                                          
\put(50,20){\line(1,1){10}}                                                 
\put(50,20){\circle*{0.8}}                                                  
\put(46,21){\scriptsize$\nabla_{2}$}                                        
\thinlines                                                                  
\put(70,20){\line(-1,1){10}}                                                
\put(70,20){\circle*{0.8}}                                                  
\put(71,21){\scriptsize$\nabla_3$}                                          
\thicklines                                                                 
\put(100,20){\line(0,1){10}}                                                
\put(100,20){\circle*{0.8}}                                                 
\put(101,21){\scriptsize$\nabla_{1}$}                                       
\put(130,21){\small$f(2)=\;5$}                                              
\thinlines                                                                  
\put(5,10){\line(1,2){5}}                                                   
\put(5,10){\circle*{0.8}}                                                   
\put(1,11){\scriptsize$\nabla_{2}$}                                         
\put(15,10){\line(-1,2){5}}                                                 
\put(15,10){\circle*{0.8}}                                                  
\put(16,11){\scriptsize$\nabla_3$}                                          
\thicklines                                                                 
\put(30,10){\line(0,1){10}}                                                 
\put(30,10){\circle*{0.8}}                                                  
\put(31,11){\scriptsize$\nabla_{1}$}                                        
\put(45,10){\line(1,2){5}}                                                  
\put(45,10){\circle*{0.8}}                                                  
\put(41,11){\scriptsize$\nabla_{2}$}                                        
\thinlines                                                                  
\put(55,10){\line(-1,2){5}}                                                 
\put(55,10){\circle*{0.8}}                                                  
\put(56,11){\scriptsize$\nabla_3$}                                          
\put(70,10){\line(0,1){10}}                                                 
\put(70,10){\circle*{0.8}}                                                  
\put(71,11){\scriptsize$\nabla_{1}$}                                        
\put(95,10){\line(1,2){5}}                                                  
\put(95,10){\circle*{0.8}}                                                  
\put(91,11){\scriptsize$\nabla_{2}$}                                        
\thicklines                                                                 
\put(105,10){\line(-1,2){5}}                                                
\put(105,10){\circle*{0.8}}                                                 
\put(106,11){\scriptsize$\nabla_3$}                                         
\put(130,11){\small$f(3)=\;8$}                                              
\thinlines                                                                  
\put(2,5){\line(3,5)3}                                                      
\put(8,5){\line(-3,5)3}                                                     
\put(15,5){\line(0,1){5}}                                                   
\put(27,5){\line(3,5)3}                                                     
\thicklines                                                                 
\put(33,5){\line(-3,5)3}                                                    
\put(42,5){\line(3,5)3}                                                     
\thinlines                                                                  
\put(48,5){\line(-3,5)3}                                                    
\put(55,5){\line(0,1){5}}                                                   
\put(67,5){\line(3,5)3}                                                     
\put(73,5){\line(-3,5)3}                                                    
\put(92,5){\line(3,5)3}                                                     
\put(92,5){\circle*{0.8}}                                                   
\put(87,6){\scriptsize$\nabla_{2}$}                                         
\put(98,5){\line(-3,5)3}                                                    
\put(98,5){\circle*{0.8}}                                                   
\put(99,6){\scriptsize$\nabla_3$}                                           
\thicklines                                                                 
\put(105,5){\line(0,1){5}}                                                  
\thinlines                                                                  
\put(105,5){\circle*{0.8}}                                                  
\put(106,6){\scriptsize$\nabla_{1}$}                                        
\put(130,6){\small$f(4)=13$}                                                
\put(89,0){\line(3,5)3}                                                     
\put(95,0){\line(-3,5)3}                                                    
\put(98,0){\line(0,1){5}}                                                   
\put(102,0){\line(3,5)3}                                                    
\thicklines                                                                 
\put(108,0){\line(-3,5)3}                                                   
\put(65,-1){\small Fig. $1$}                                                
\put(130,1){\small$f(5)=21$}                                                
\thinlines                                                                  
\end{picture} 

\medskip
\noindent
Let $f_{i}(k)$ be a number of meaningful compositions of the
$k^{\mbox{\scriptsize th}}$-order beginning with $\nabla_{i}$.
Let $f(k)$ be a number of meaningful composition of the
$k^{\mbox{\scriptsize th}}$-order of operations over $\Theta$.
Then $f(k)=f_1(k) + f_2(k) + f_3(k)$. Based on partial self similarity
of the tree~(Fig.~$\!1$), which is formed according to {\sc Cayley}'s
table (4), we get equalities:
{\small $$
f_{1}(k) = f_{2}(k-1) + f_3(k-1) \;\; \wedge \;\;
f_{2}(k) = f_{2}(k-1) + f_3(k-1) \;\; \wedge \;\;
f_3(k) = f_{1}(k-1).
$$}
\hspace*{-2.5 mm} Now, a recurrent relation for $f(k)$
can be derived as follows:
{\small $$
\begin{array}{rcl}
f(k) &\!\!\!=\!\!\!& f_{1}(k) + f_{2}(k) + f_3(k)                      \\
     &\!\!\!=\!\!\!& {\big (} f_{1}(k-1) + f_{2}(k-1) + f_3(k-1) {\big )}
            + {\big (} f_3(k-1) + f_{2}(k-1) {\big )}                  \\
     &\!\!\!=\!\!\!& f(k-1) + {\big (} f_{1}(k-2) + f_{2}(k-2) + f_3(k-2) {\big )}
            = f(k-1) + f(k-2).
\end{array}
$$}
\hspace*{-2.5 mm} Based on the initial values: $f(1)=3, \, f(2)=5, \, f(3)=8$
we conclude that $f(k) = F_{k+3}$, where is {\sc Fibonacci}'s number of
order $k+3$.

\medskip
Let us note that $\nabla_2\circ\nabla_1 = 0$ and $\nabla_3\circ\nabla_2 = 0$,
because $d^2 = 0$. On the other hand, the compositions $\nabla_1\circ\nabla_3$,
$\nabla_2\circ\nabla_2$ and $\nabla_3\circ\nabla_1$ are not annihilated,
because of $\varphi_{0}^{-1}\circ\varphi_{3} \neq i$ and $\varphi_{1}^{-1}
\circ \varphi_{2} \neq i$. Thus, as in the paper [{\bf 2}], we conclude that
the non-trivial compositions are of the following form:
$$                                                                       
\begin{array}{c}
( \nabla_1 \circ ) \nabla_3 \circ \cdots \circ \nabla_1 \circ \nabla_3 \circ \nabla_1, \\
  \nabla_2 \circ   \nabla_2 \circ \cdots \circ \nabla_2 \circ \nabla_2 \circ \nabla_2, \\
( \nabla_3 \circ ) \nabla_1 \circ \cdots \circ \nabla_3 \circ \nabla_1 \circ \nabla_3.
\end{array}
\leqno (5)
$$
As non-trivial compositions we consider those which are not identical to the
zero function. Terms in parentheses are included in for an odd number of
terms and are left out otherwise.

\break

\bigskip

\centerline{\bf 2. DIFFERENTIAL FORMS AND OPERATIONS ON THE SPACE
${\mbox{\bf R}}^{\mbox{\footnotesize \bf n}}$}

\bigskip

Let us present a recurrent relation for counting meaningful compositions of
the higher-order differential operations on the space $\mR^n$ $(n=3,4,\ldots)$
and extract the non-trivial compositions of order higher than two. Let
us form the following sets of functions:
{\small $$
\mbox{\normalsize $\mA_{i}$} = \{ \mbox{\boldmath $f$} : \mR^n \rightarrow
\mR^{n \choose i} | f_1, \ldots, f_{n \choose i} \in C^{\infty}(\mR^{n}) \}
$$}
\hspace*{-2.5 mm} for $i = 0, 1, \ldots, m$ where $m = [n/2]$. Let $\Omega^i
(\mR^n)$ be a set of differential forms of degree $i = 0, 1, \ldots, n$ on
the space $\mR^n$. Let us notice that $\Omega^i(\mR^n)$ and $\Omega^{n-i}(\mR^n)$,
over ring $\mA_{0}$, are spaces of the same dimension $n \choose i$, for $i = 0, 1,
\ldots, m$. They can be identified with $\mA_i$, using the corresponding
isomorphisms:
$$
\varphi_{i} :
\Omega^{i}(\mR^{n}) \rightarrow \mA_{i}
\;\; (\mbox{\small $0$} \leq i \leq \mbox{\small $m$})
\;\;\;\; \mbox{and} \;\;\;\;
\varphi_{n-i} :
\Omega^{n-i}(\mR^{n}) \rightarrow \mA_{i}
\;\; (\mbox{\small $0$} \leq i < \mbox{\small $n-m$}).
$$

\noindent
\begin{minipage}[b]{22 cc}
We define the first-order differential operations on the space $\mR^n$
via the operator of the exterior differentiation $d$ as follows:
$$
\hspace*{20 mm} \nabla_{i} = \varphi_{i} \circ d \circ \varphi_{i-1}^{-1}
\;\; (\mbox{\small $1$} \leq i \leq \mbox{\small $n$}).
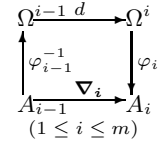
$$
\end{minipage}
{\small \unitlength 0.76 mm 
\begin{picture}(16,16)(-3,-3)                                                
\put(2,-3){\scriptsize $(1 \leq i \leq m)$}                                  
\put(0,16){$\Omega^{i-1}$}                                                   
\put(1,4){\vector(0,1){11}}                                                  
\put(2,9){\scriptsize $\varphi^{-1}_{i-1}$}                                  
\put(0,1){$A_{i-1}$}                                                         
\put(3,3){\vector(1,0){16}}                                                  
\put(10,4){\scriptsize \boldmath $\nabla_i$}                                 
\put(19,1){$A_{i}$}                                                          
\put(20,15){\vector(0,-1){11}}                                               
\put(21,9){\scriptsize $\varphi_{i}$}                                        
\put(19,16){$\Omega^{i}$}                                                    
\put(3,17){\vector(1,0){16}}                                                 
\put(10,18){\scriptsize $d$}                                                 
\end{picture} } 

\smallskip
\noindent
Therefore, we obtain the first order differential operations on the space
$\mR^n$, depending on pairity of dimension $n$, in the following form:
\begin{center}
\footnotesize
\begin{tabular}{cc}
$ \begin{array}{ll}
\mbox{\small $n=2m:$}
           & \mbox{\small $\nabla_{1}$}   : \mA_{0} \rightarrow \mA_{1}   \\
           & \mbox{\small $\nabla_{2}$}   : \mA_{1} \rightarrow \mA_{2}   \\
           & \vdots                                                       \\
           & \mbox{\small $\nabla_{i}$}   : \mA_{i} \rightarrow \mA_{i+1} \\
           & \vdots                                                       \\
           & \mbox{\small $\nabla_{m}$}   : \mA_{m-1} \rightarrow \mA_{m} \\
           & \mbox{\small $\nabla_{m+1}$} : \mA_{m} \rightarrow \mA_{m-1} \\
           & \vdots                                                       \\
           & \mbox{\small $\nabla_{n-j}$} : \mA_{j+1} \rightarrow \mA_{j} \\
           & \vdots                                                       \\
           & \mbox{\small $\nabla_{n-1}$} : \mA_{2} \rightarrow \mA_{1}   \\
           & \mbox{\small $\nabla_{n}$}   : \mA_{1} \rightarrow \mA_{0}
             \mbox{\normalsize ,}
\end{array} $ & $ \begin{array}{ll}
\mbox{\small $n=2m+1:$}
           & \mbox{\small $\nabla_{1}$}   : \mA_{0} \rightarrow \mA_{1}   \\
           & \mbox{\small $\nabla_{2}$}   : \mA_{1} \rightarrow \mA_{2}   \\
           & \vdots                                                       \\
           & \mbox{\small $\nabla_{i}$}   : \mA_{i} \rightarrow \mA_{i+1} \\
           & \vdots                                                       \\
           & \mbox{\small $\nabla_{m}$}   : \mA_{m-1} \rightarrow \mA_{m} \\
           & \mbox{\small $\nabla_{m+1}$} : \mA_{m} \rightarrow \mA_{m}   \\
           & \mbox{\small $\nabla_{m+2}$} : \mA_{m} \rightarrow \mA_{m-1} \\
           & \vdots                                                       \\
           & \mbox{\small $\nabla_{n-j}$} : \mA_{j+1} \rightarrow \mA_{j} \\
           & \vdots                                                       \\
           & \mbox{\small $\nabla_{n-1}$} : \mA_{2} \rightarrow \mA_{1}   \\
           & \mbox{\small $\nabla_{n}$}   : \mA_{1} \rightarrow \mA_{0}
             \mbox{\normalsize .}
\end{array} $
\end{tabular}
\end{center}
Consider the set of functions $\Theta = \{ \nabla_1, \nabla_2, \ldots,
\nabla_n \}$. Let us define a binary relation $\rho$~{\it "to be in
composition"} with $\nabla_i \rho \nabla_j = \top$ iff the composition
$\nabla_{j} \circ \nabla_{i}$ is meaningful $(\nabla_i, \nabla_j \in \Theta)$.
It is not difficult to check that {\sc Cayley}'s table of this relation
is determined with:
{\small
$$                                                                      
\mbox{\normalsize $\nabla_{i} \rho \nabla_{j}$} =
\left\{
\begin{array}{lll}
\top &:& (j = i + 1)    \vee   (i + j = n + 1),  \\
\bot &:& (j \neq i + 1) \wedge (i + j \neq n + 1).
\end{array}
\right.
\leqno \mbox{\normalsize (6)}
$$}
\hspace*{-2.5 mm}
Let us form an adjacency matrix $\mbox{\tt A}=[a_{ij}]\in\{0,1\}^{n\times n}$
of the graph, determined by relation $\rho$. Let $f_{i}(k)$ be a number of
meaningful compositions of the $k^{\mbox{\scriptsize th}}$-order beginning with
$\nabla_{i}$ (notice that $f_{i}(1)\!=\!1$ for $i\!=\!1,\ldots,n$). Let $f(k)$
be a number of meaningful composition of the $k^{\mbox{\scriptsize th}}$-order
of operations over $\Theta$. Then $f(k)\!=\!f_{1}(k)\!+\!\ldots\!+\!f_{n}(k)$.
Notice that the following is true:
$$                                                                      
f_{i}(k) = \displaystyle \sum_{j=1}^{n}{a_{ij} \cdot f_{j}(k-1)},
\leqno (7)
$$
for $i = 1,\ldots,n$. Based on $(7)$ we form the system of recurrent equations:
$$                                                                      
\left[
\begin{array}{c}
f_{1}(k) \\
\vdots   \\
f_{n}(k)
\end{array}
\right]
=
\left[
\begin{array}{ccc}
a_{11} & \cdots & a_{1n} \\
\vdots &        & \vdots \\
a_{n1} & \cdots & a_{nn}
\end{array}
\right]
\cdot
\left[
\begin{array}{c}
f_{1}(k-1) \\
\vdots     \\
f_{n}(k-1)
\end{array}
\right]_{\mbox{\normalsize .}}
\leqno (8)
$$
If $v_{n} = [ \; 1 \;\; \cdots \;\; 1 \; ]_{1 \times n}$ then:
{\small $$                                                              
\mbox{\normalsize $f(k)$}
\;\mbox{\normalsize $=$}\;
\mbox{\normalsize $v_{n}$}
\cdot
\left[
\begin{array}{c}
\mbox{\normalsize $f_{1}(k)$} \\
\vdots                        \\
\mbox{\normalsize $f_{n}(k)$}
\end{array}
\right]_{\mbox{\normalsize .}}
\leqno \mbox{\normalsize $(9)$}
$$}
\hspace*{-2.5 mm} So, the expression:
$$                                                                     
f(k) = v_{n} \cdot \mbox{\tt A}^{k-1} \cdot v^{T}_{n}.
\leqno (10)
$$
follows from $(8)$ and $(9)$. Reducing the system of the recurrent equations
$(8)$, for any of the functions $f_{i}(k)$ we have:
$$                                                                     
\alpha_{0}f_{i}(k)+\alpha_{1}f_{i}(k-1)+\cdots+\alpha_{n}f_{i}(k-n)=0
\qquad (k > n),
\leqno (11)
$$
where $\alpha_{0}, \ldots, \alpha_{n}$ are coefficients of the characteristic
polynomial $P_{n}(\lambda) = |\mbox{\tt A} - \lambda \mbox{\tt I}| = \alpha_{0}
\lambda^{n} + \ldots + \alpha_{n}$. Thus, we conclude that the function $f(k) =
\sum\limits_{i=1}^{n}{f_{i}(k)}$ also satisfies:
$$                                                                     
\alpha_{0}f(k)+\alpha_{1}f(k-1)+\cdots+\alpha_{n}f(k-n)=0
\qquad (k > n).
\leqno (12)
$$

\noindent
Hence, the following theorem holds.

\bigskip
\noindent
{\bf Theorem 1.} {\it The number of meaningful differential operations,
on the space $\mR^{n}$ $(n=3,4,\ldots)$, of the order higher than two,
is determined by the formula $(10)$, i.e. by the recurrent formula $(12)$.}

\bigskip
In $n$-dimensional space $\mR^n$, for dimensions $n = 3, 4, 5, \ldots, 10$,
using the previous theorem we form a table of the corresponding recurrent
formula:

{\footnotesize
\begin{center}
\begin{tabular}{|c|c|}
\hline {\small \rm Dimension}:
 &     {\small \rm \quad Recurrent relations for
 the number of meaningful compositions: \quad}                       \\ \hline
$n = \;$ 3 & $f(i+2)=f(i+1)+f(i)$                                    \\ \hline
$n = \;$ 4 & $f(i+2)=2 f(i)$                                         \\ \hline
$n = \;$ 5 & $f(i+3)=f(i+2) + 2 f(i+1) - f(i)$                       \\ \hline
$n = \;$ 6 & $f(i+4)=3 f(i+2) - f(i)$                                \\ \hline
$n = \;$ 7 & $f(i+5)=f(i+3) + 3 f(i+2) - 2  f(i+1) - f(i)$           \\ \hline
$n = \;$ 8 & $f(i+4)=4 f(i+3) - 3 f(i)$                              \\ \hline
$n = \;$ 9 & $f(i+5)=f(i+4) + 4 f(i+3) - 3 f(i+2) - 3 f(i+1) + f(i)$ \\ \hline
$n = 10$ & $f(i+6)=5 f(i+4) - 6 f(i+2) + f(i)$                       \\ \hline
\end{tabular}
\end{center}}

\break

\smallskip
Let us determine non-trivial higher-order meaningful compositions on the
space $\mR^n$. For isomorphisms $\varphi_k$ we have:
$$                                                                      
\varphi_{k}^{-1} \circ \varphi_{n-k} \neq i,
\leqno (13) $$
for $k = 1 , 2, \ldots, n$ and  $2k \neq n$. Then, based on (6) and (13),
all second-order compositions are given by the formula:
$$                                                                      
\nabla_{j} \circ \nabla_{k} =
\left\{ \begin{array}{lll}
0         \!\!&:&\!\! j = k + 1,                              \\
g_{j,k}   \!\!&:&\!\! (k + j = n + 1) \wedge (2 k \neq n),    \\
\vartheta \!\!&:&\!\! (j \neq k + 1)  \wedge (k + j \neq n + 1);
\end{array}
\right.
\leqno (14) $$
where $0$ is a trivial composition, $g_{j,k}$ is a non-trivial second-order
composition and $\vartheta$ is a nowhere-defined function for $j, k = 1,
\ldots, n$. Notice that in $g_{j,k} = \nabla_{j} \circ \nabla_{k} =
\varphi_{n+1-k}\!\:\circ\!\:d\!\:\circ\!\:\varphi_{n-k}^{-1}\!\:\circ\!\:
\varphi_{k}\!\:\circ\!\:d\!\:\circ\!\:\varphi_{k-1}^{-1}$
$(j\!=\!n\!+\!1\!-\!k \: \wedge \: 2 k\!\neq\!n)$ and switching
the terms is impossible, because in that way we get nowhere-defined
function $\vartheta$. Hence, we conclude that the following theorem
holds.

\bigskip
\noindent
{\bf Theorem 2.} {\it All meaningful non-trivial differential operations on
the space $\mR^n$ $(n=3,4,\ldots)$, of order higher than, two are given in
the~form~of~the~following~compositions}:
$$                                                                      
\begin{array}{c}
(\nabla_{k}) \circ \nabla_j \circ \nabla_k \circ \cdots \circ \nabla_j \circ \nabla_k, \\
(\nabla_{j}) \circ \nabla_k \circ \nabla_j \circ \cdots \circ \nabla_k \circ \nabla_j,
\end{array}
\leqno (15) $$
{\it with to the condition $k + j = n + 1$ and $2 k, \: 2j \neq n$ for
$k, \: j = 1, 2, \ldots, n$. Terms in parentheses are included in for an
odd number of terms and are left out otherwise.}

\bigskip
\noindent
{\bf Acknowledgment.} I wish to express my gratitude to Professors
{\sc M. Merkle} and {\sc M. Prvanovi\' c} who examined the first version
of the paper and gave me their suggestions and some very useful remarks.

\bigskip

\bigskip

\bigskip
\begin{center} \small \bf REFERENCES \end{center}

\bigskip

\newcounter{ref} \begin{list}{\small \arabic{ref}.}
{\usecounter{ref} \leftmargin 4mm \itemsep -1mm}

\item {\small {\sc R.$\:$Bott, L.$\:$W.$\:$Tu:} {\it Differential forms
in algebraic topology}, Springer, New York 1982.}

\item {\small {\sc B.$\:$J.$\:$Male\v sevi\' c:} {\it A note on higher-order
differential operations}, Univ. Beograd, Publ. Elektrotehn. Fak.,Ser. Mat.
{\bf 7} (1996), 105-109.}
\end{list}

\bigskip

\bigskip
{\small
\noindent University of Belgrade,
           \hfill (Received September 8, 1997)       \break
\noindent Faculty of Electrical Engineering,
           \hfill (Revised October 30, 1998)         \break
\noindent P.O.Box 35-54, $11120$ Belgrade,     \hfill\break
\noindent Yugoslavia                                 \break
\noindent {\footnotesize \bf malesevic@kiklop.etf.bg.ac.yu}
\hfill}

\newpage

\end{document}